\documentclass[12pt]{article}
\usepackage{graphics}
\date{}
\def\qed{{\hfill\rule{1.2ex}{1.2ex}}}
\newtheorem{Th}{Theorem}

\def\N{{\bf N}}

\def\Tr{\mathop{\hbox{\rm Tr}}\nolimits}

\author{Norbert A'Campo}
\title{A combinatorial property of generic immersions of curves}

\begin{document}

\maketitle

A {\it divide} (called ``partage'' in [AC2]) with $r$ branches is the image $P$ of a
generic relative immersion of
the union of $r, r \in \N,$
copies of the unit interval $[0,1]$ in the unit disk $D^2$. A {\it region} of
a divide $P$ is a connected component $A$ of $D^2 \setminus P$,
that does not meet the boundary $\partial{D^2}$. A {\it segment} of $P$ is a 
connected component of the complement of the double points of $P$ in $P$. 
A {\it sector} of $P$ is a connected component of the germ 
at a double
point of a region. We say that
a divide is {\it cellular}, if $P$ is connected and the closure in $D^2$ of
each region of $P$ is contractible. The {\it link}
$L(P) \subset S^3$ of a divide $P$ is obtained by a hodographic 
construction [AC4-5]:
$$
L(P):=\{(x,u) \in TD^2 \mid x \in P, u \in T_xP, ||x||^2+||u||^2=1 \}
$$
The link $L(P)$ of a divide has a natural orientation.
The complement $S^3 \setminus L(P)$ of the link of a connected divide $P$
admits a fibration over the circle $S^1$, whose restriction to each oriented
meridian of each component is a degree $1$ map [AC5]. The fiber $F_P$ 
is a surface
of genus $\delta(P)$ with $r$ boundary components, where $\delta(P)$ is
the number of double points of $P$. Let $T_P:F_P \to F_P$ be the monodromy
diffeomorphism of this fibration. A divide $P$ is {\it simple} 
if $P$ is connected
and has at least one double point, such that 
there does not exist a relatively embedded copy of $[0,1]$ in $D^2$, that
cuts $P$ transversally in one point and separates the double points
of $P$ in two non-empty sets.

\begin{Th} The Lefschetz number of the monodromy $T_P$ of a simple, 
cellular divide is $0$.
\end{Th}

If the cellular divide $P$ is the saddle level of a local real 
maximal deformation 
of a real equation of a  plane curve
singularity [S1-2,AC2,AC4,G-Z], the Lefschetz number of the 
monodromy $T_P$ equals the
Lefschetz number of the monodromy of the singularity, and therefore 
equals $0$ by [AC1].
The divide of Fig. $1$ is obtained by a maximal real deformation 
from the singularity
$x((y^2-x^3)^2-4yx^4-x^7)$ [G-Z]. So the Lefschetz number of 
its monodromy equals
$0$ although the divide  is not cellular.

\begin{center}
\scalebox{1}[1]{\includegraphics{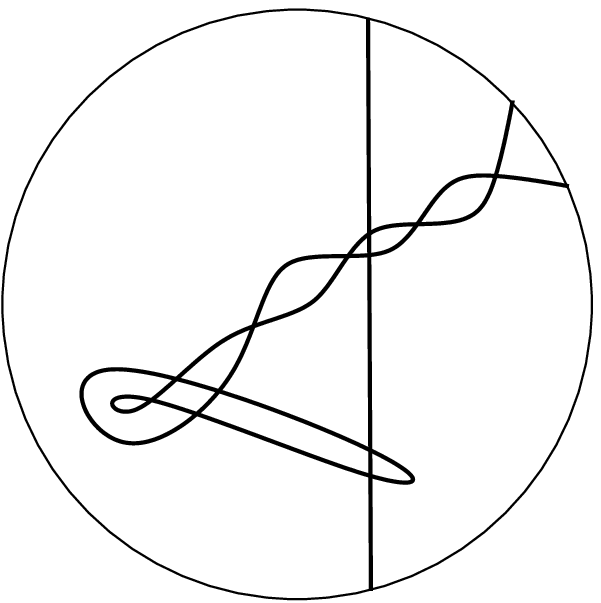}}
\newline
{Fig. 1. Divide for the singularity $x((y^2-x^3)^2-4yx^4-x^7)$.}
\end{center}

The Lefschetz
number of the non-cellular divide of Fig. $2a$ equals $2$. The Lefschetz
number of the non-simple divide of Fig. $2b$ equals $-1$.

\begin{center}
\scalebox{1}[1]{\includegraphics{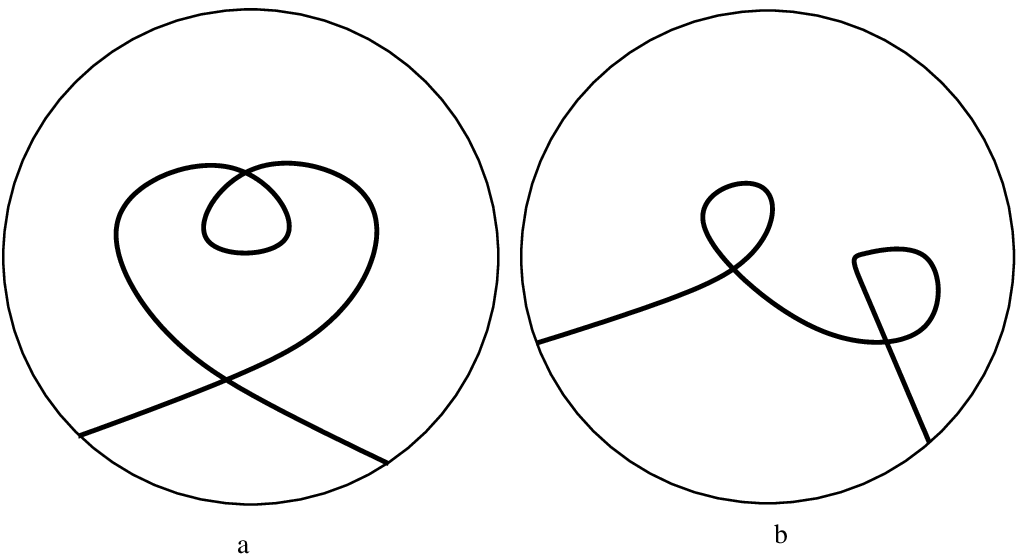}}
\newline
{Fig. 2. a: non-cellular, b: non-simple.}
\end{center}

Before giving the proof, we will recall some material about the Seifert
form of the link of a connected divide, see [AC2,AC4, G-Z]. The connected 
components of $D^2 \setminus P$
will be signed with $\pm$, such that two components have opposite signs
if their closures have a segment of $P$ in common. Moreover, we choose 
a basepoint $b_A$ in 
each region $A$. The {\it geometric Dynkin 
diagram} $\Gamma$ of the divide $P$ is the
planar graph in $D^2$ obtained as follows. The vertices of $\Gamma$ are 
the double points of $P$ and the basepoints of the regions.
The edges of $\Gamma$ are arcs of two species. For each sector 
at a double point of $P$ we
draw  an arc, that connects the double point with
the basepoint of the region to which the sector belongs. Moreover, for 
each segment
of $P$, we connect the basepoints of those regions, whose 
closures meet along the given segment by an arc having one 
transversal intersection with the segment. 
Pairs of vertices of $\Gamma$ can be connected  by several edges precisely
if the divide is non-cellular. 
	
The geometric Dynkin diagram was 
introduced by Sabir Gusein-Zade [G-Z] as system of gradient lines of a morse 
function on $D^2$.
We number the vertices of 
$\Gamma$ with $1\leq i \leq \mu_P$ by taking first the basepoints of the
$-$ regions, then the double points of $P$ and finally the basepoints of the
$+$ regions. The geometric Dynkin diagram  defines the Seifert
form $S$ of the fibered link $L(P)$. On the vector space spanned by the
set of vertices $V$ of $\Gamma$ we define an upper triangular 
nilpotent endomorphism
$N:=(n_{ij})$ by setting  the matrix coefficient $n_{ij}$  
for $i<j$ equal
to the number of edges of $\Gamma$ that connect $i$ and $j$.

From the tricollaring of $\Gamma$ by $\{-,.,+\}$ one can conclude 
for any divide a basic fact, namely  
that the equality $N^3=0$ holds. A 
divide is a slalom divide of a rooted tree or of a disk-wide-web 
if and only if the equality $N^2=0$ holds, see [AC3,AC6].  

The Seifert form is represented by $Id+N$. 
The monodromy action $T_*$ on the first
homology of $F_P$ is represented by the matrix $^t{(Id+N)}^{-1}\circ (Id+N)$.
For the Lefschetz number $\Lambda_P$ of $T$ this yields 
the expression
$$
1-\Tr(Id-^t{N}+N-^t{N}N+^t{N}^2N)=1-\mu_P+\Tr(^t{N}N)-\Tr(^t{N}^2N)
$$
For instance for a slalom divide we have
$$
\Lambda=1-\mu_P+\Tr(^t{N}N)=1-\mu_P+(\mu_P-1)=0.
$$

\noindent
{\bf Proof of  theorem $1$:}
For a cellular divide we have $n_{ij}^2=n_{ij}$, hence 
$$
\Tr(^t{N}N)=\sum_{ij}n_{ij}^2=\sum_{ij}n_{ij}=e_P
$$ 
where
$e_P$ is the number of edges of the geometric Dynkin diagram $\Gamma$. The term 
$\Tr(^t{N}^2N)$ has for a cellular divide also a 
combinatorial interpretation. We call flag 
of the divide $P$ a pair $(a,b)$ of edges of $\Gamma$ such that the
edge $a$ connects a basepoint of a $-$ region with a 
double point of $P$, while the edge $b$ connects this double point 
with the basepoint of a $+$ region. Calling the number of flags $f_P$
we have for a cellular divide $\Tr(^t{N}^2N)=f_P$ and it remains to prove:
$\mu_P-e_P+f_P=1.$

\begin{center}
\scalebox{1}[1]{\includegraphics{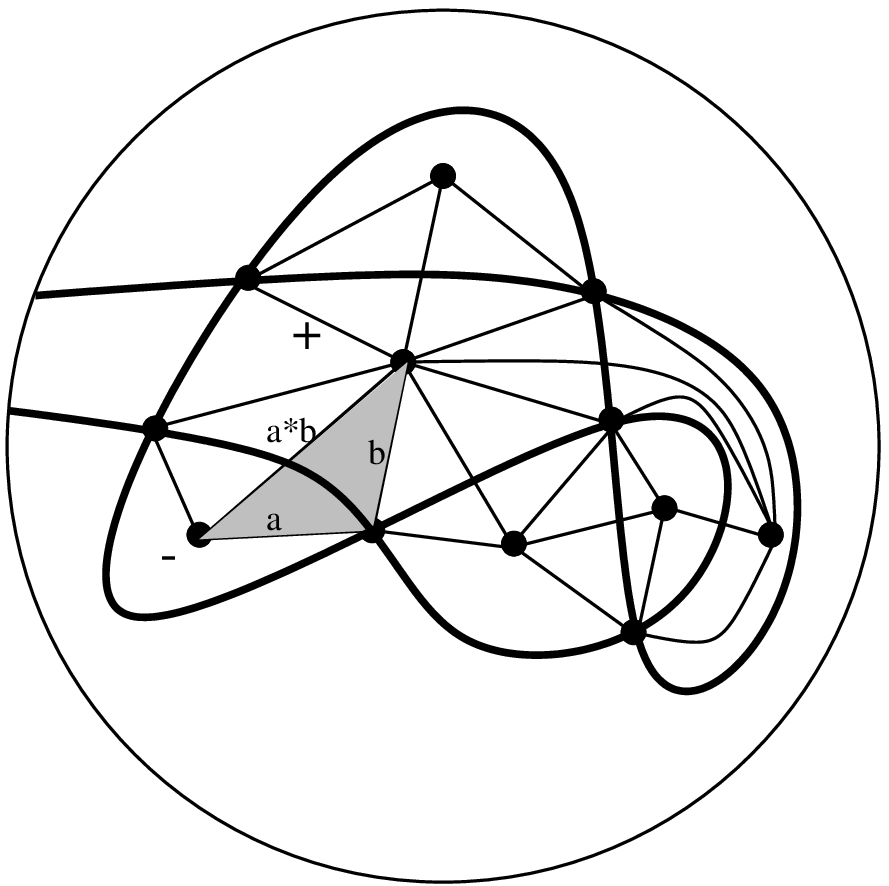}}
\newline
{Fig. $3$. Body of a divide with its partial triangulation.}
\end{center}

We call {\it body} $B_P$ of the divide $P$ 
the union of the set of its double points with
the closure of the union of its regions. For a simple 
divide $P$ the Euler-Poincar\'e
characteristic $\chi(B_P)$ equals $1$. Each flag $(a,b)$ of $P$ 
defines a triangle $\Delta_{(a,b)}=(a,b,a*b)\subset B_P$ 
in the the body,  where $a*b$ is the edge
of $\Gamma$ which connects the non common endpoints of $a$ and $b$. 
Let $S \subset B_P$ be the union of the set 
of double points of $P$, of the edges of $\Gamma$ and of the triangles
associated to the flags. This system of 
vertices, edges and triangles is a triangulation of  $S$. The body $B_P$ 
collapses onto $S$, see Fig $3$. Hence we conclude
$$
\mu_P-e_P+f_P=\chi(S)=\chi(B_P)=1,
$$
which completes  the proof.

I like to thank Marc Burger for his remark, that $\mu_P-e_P+f_P$
should be interpreted as Euler-Poincar\'e characteristic.

\noindent
{\bf Remark:} The property $N^3=0$ that, as we have seen here above, 
holds for any divide, indicates
that the computation of the traces of the iterates of the monodromy
is related to the random walk on the geometric Dynkin diagram $\Gamma$.
We ask for an expression of $\Tr(T^k)$ in terms of the generating
function of the random walk on the Dynkin diagram $\Gamma$, see [D-J].

\noindent
[AC1]
Norbert A'Campo,
{\it Le nombre de Lefschetz d'une monodromie},
Nederl. Akad. Wetensch. Proc. Ser. A \{76\}=Indag. Math. 
{\bf 35} (1973), 113--118,

\noindent
[AC2]
Norbert A'Campo,
{\it Le Groupe de Monodromie du D\'eploiement des Singularit\'es
Isol\'ees de Courbes Planes I},
Math. Ann. {\bf 213} (1975)
1--32.

\noindent
[AC3]
Norbert A'Campo,
{\it Le Groupe de Monodromie du D\'eploiement des Singularit\'es
Isol\'ees de Courbes Planes II}, 
Actes du Con\-gr\`es Inter\-national des Math\'ema\-ti\-ciens,  
tome 1,
395--404,
Vancouver,
1974.

\noindent
[AC4]
Norbert A'Campo,
{\it Real deformations and complex topology of plane curve singularities},
Annales de la Facult\'e des Sciences de Toulouse {\bf 8} (1999), 1, 5--23.

\noindent
[AC5]
Norbert A'Campo,
{\it Generic immersions of curves, knots,
monodromy and gordian number},
Publ. Math. I.H.E.S. {\bf 88} (1998), 151-169, (1999).

\noindent
[AC6]
Norbert A'Campo,
{\it Planar trees, slalom curves and hyperbolic knots},
Publ. Math. I.H.E.S. {\bf 88} (1998), 171-180, (1999).

\noindent
[D-J]
E.B. Dynkin, A.A. Juschkewitsch, {\it Saetze und Aufgaben ueber Markoffsche Prozesse}, aus dem Russischen uebersetzt von Klaus Sch\"urger: {\it Teoremy i zadachi o prochessakh Markova}, Moskva 1967,  Berlin : Springer, 1969 
Heidelberger Taschenb"ucher.

\noindent
[G-Z]
S. M. Gusein-Zade,
{\it Matrices d'intersections pour certaines singularit\'es de
fonctions de 2 variables},
Funkcional. Anal. i Prilozen
{\bf 8}, (1974),
11--15.

\noindent
[S1]
Charlotte Angas Scott,
{\it On the Higher Singularities of Plane Curves},
Amer. J. Math. {\bf 14},(1892) 301--325.

\noindent
[S2]
Charlotte Angas Scott,
{\it The Nature and Effect of Singularities of Plane Algebraic Curves},
Amer. J. Math. {\bf 15}, (1893) 221--243.

\noindent
Universitaet Basel, Rheinsprung 21, CH-4051  Basel.
\end{document}